\documentclass[12pt]{amsart}

\newcommand{\RP}{{\bf RP}}
\newcommand{\Z}{{\bf Z}}
\newcommand{\R}{{\bf R}}

\newcommand{\e}{\varepsilon}
\newcommand{\x}{{\bf x}}
\newcommand{\sm}{\setminus}
\newcommand{\F}{{\mathcal F}}

\begin{document}

\author{V.~A.~Vassiliev}

\email{vva@mi.ras.ru}

\thanks{
Supported by RFBR (project 95-01-00846a), INTAS (project 4373) and Netherlands
Organization for Scientific Research (NWO), project 47.03.005.}

\title{On $r$-neighborly submanifolds in $\R^N$}

\date{Revised version was published in 1998 and dedicated to Professor
J\"urgen Moser on the occasion of his 70-th birthday}

%\begin{flushright}
%\centerline{To }
%\end{flushright}

\begin{abstract}
A submanifold $M \subset \R^N$ is $r$-{\it neighborly} if  for any $r$ points
in $M$ there is a hyperplane, supporting $M$ and touching it at exactly these
$r$ points. We prove that the minimal dimension $\Delta(k,r)$ of the Euclidean
space, containing a stably $r$-neighborly submanifold, is asymptotically no
less than $2kr-k$.  \end{abstract}

\maketitle

\section{Introduction}

Let $M$ be a $k$-dimensional manifold, and $r$ a natural number.

\medskip
{\sc Definition.} A smooth embedding $M \to \R^N$ is $r$-{\it neighborly} if
for any $r$ points in $M$ there is an affine hyperplane in $\R^N$, supporting
$M$ and touching it at exactly these $r$ points.
\medskip

Denote by $\delta(M,r)$ the minimal dimension $N$ of an Euclidean space, such
that there exists an $r$-neighborly embedding $M \to \R^N$, and by
$\delta(k,r)$ the maximum of numbers $\delta(M,r)$ over all connected
$k$-dimensional manifolds $M$.
\medskip

The problem of estimating the numbers $\delta(k,r)$ for all $k$ and $r$ was
posed by M.~Perles in the 1970-ies by analogy with similar problems of
combinatorics, and was discussed in the Oberwolfach Combinatorics meetings in
1980-ies.
%, see \cite{prlist}.
Nontrivial examples of $r$-neighborly manifolds were constructed in \cite{HM};
as G.~Kalai communicated to me, in their non-published work with A.~Wigderson a
polynomial upper estimate of type $\delta(k,r)=O(r^2 k)$ was proved.  However
no nontrivial general lower estimates of these numbers are known.
\medskip

We consider the similar problem concerning a slightly stronger condition.

\medskip
{\sc Definition.} A smooth embedding $M \to \R^N$ is {\it stably} $r$-{\it
neighborly}, if it is $r$-neighborly, and any other embedding, sufficiently
close to it in the $C^2$-topology, also is. The corresponding analogues of
numbers $\delta(M,r)$ and $\delta(k,r)$ are denoted by $\Delta(M,r)$ and
$\Delta(k,r)$, respectively.
\medskip

It is more or less obvious that
\begin{equation}
\label{obv} \Delta(k,r) \ge (k+1)r.
\end{equation}

Indeed, if $N<(k+1)r,$ then for any generic submanifold $M^k \subset \R^N$ and
almost any set of r points in $M^k$  the minimal affine plane in $\R^N,$
touching $M^k$ at all points of this set, coincides with entire $\R^N.$

If $k=1$, then this estimate is sharp: the moment embedding $S^1 \to \R^{2r}$,
sending the point $\alpha$ to $(\sin \alpha, \cos \alpha, \ldots, \sin r\alpha,
\cos r\alpha),$ is stably $r$-neighborly.
\medskip

We prove the following lower estimates of numbers $\Delta(M,r)$. Denote by
$d(r)$ the number of ones in the binary representation of $r$.
\medskip

{\sc Theorem 1.} {\it If $k$ is a power of 2, then}
\begin{equation}
\label{th1} \Delta(\R^k,r) > (k+1)r+(k-1)(r-d(r))-1.
\end{equation}
{\it In particular, if $r$ also is a power of 2, then} $\Delta(\R^k,r) >
2kr-k.$
\medskip

Of course, the same estimate of the number $\Delta(M,r)$ is true for any
$k$-dimensional manifold $M$. \medskip

For an {\it arbitrary fixed} $r$ and growing $k$ we have a slightly better
estimate of $\Delta(k,r)$.
\medskip

{\sc Theorem 2.} {\it For any} $r \le k = 2^j,$
\begin{equation}
\label{th2} \Delta(\RP^k,r) > (k+1)r+kr-r(r+1)/2-1.
\end{equation}

The estimate of Theorem 2 can be improved if $r=1$.
\medskip

{\sc Proposition 1.} {\it For any closed $M$, $\Delta(M,1) \le N$ if and only
if $M$ can be smoothly embedded into $S^{N-1}.$ In particular, for} $k=2^j$,
$\Delta(\RP^k,1) = 2k+1.$ \medskip

{\sc Conjecture.} In both theorems 1 and 2, $\Delta$ can be replaced by
$\delta$.
\medskip

In fact, we prove this conjecture in a much more general situation than it
follows from Theorems 1 and 2.
\medskip

{\sc Theorems 1!, 2!}. {\it All $r$-neighborly embeddings $\R^k \to \R^N$
(respectively, $\RP^k \to \R^N$) with $N$ equal to the right-hand part of the
inequality (\ref{th1}) (respectively, (\ref{th2})), belong to a subset $\Sigma
\Sigma$ of infinite codimension in the space $C^2(\R^k,\R^N)$ (respectively,
$C^2(\RP^k,\R^N)$).}
\medskip

I believe that in fact such embeddings do not exist. \medskip

Of course, Theorems 1 and 2 follow from these ones and even from similar
statements with ``infinite codimension'' replaced by ``positive codimension''.
The set $\Sigma \Sigma$ will be described in \S \ 3.
\medskip

I thank Professor Gil Kalai for stating the problem.

\section{A topological lemma}

{\sc Convention and notation.} Using the parallel translations in $\R^N$, we
will identify the space $\R^N$ with its tangent spaces at all points $y \in
\R^N$, and use the canonical projection $\pi: T\R^N \to \R^N,$ sending any
vector $v \in T_y \R^N$ to the point $y+v$. In particular, any vector subspace
$\Xi \subset T_y \R^N$ defines the affine subspace $\pi \Xi \subset\R^N$ of the
same dimension. The sign $\oplus$ denotes the direct (Whitney) sum of two
vector bundles. $w(\xi)$ is the total Stiefel--Whitney class of the vector
bundle $\xi$, see \cite{MS}, and $\bar w(\xi)$ is the total Stiefel--Whitney
class of any vector bundle $\eta$ such that the bundle $\xi \oplus \eta$ is
trivial. In particular, $w(\xi) \bar w(\xi) = {\bf 1} \in H^0($the base of
$\xi, \Z_2)$.
\medskip

Let $L$ be a $l$-dimensional compact manifold without boundary (maybe not
connected), and $\xi$ a $t$-dimensional vector bundle over $L$ ($t < N-l$).
\medskip

{\sc Definition.} A {\it smooth homomorphism} $i: \xi \to T\R^N$ is a pair,
consisting of a smooth map $\iota:L \to \R^N$ and, for any $x \in L$, a
homomorphism $\xi_x \to T_{\iota (x)}\R^N$ smoothly depending on $x$. $i$ is a
{\it monomorphism} if all these homomorphisms $\xi_x \to T_{\iota(x)}\R^N$ are
injective.
\medskip

Let $i: \xi \to T\R^N$ be a smooth monomorphism such that the corresponding map
$\iota:L \to \R^N$ is a smooth embedding and for any $x \in L$ two subspaces
$i(\xi_x)$ and $T_{\iota(x)}(\iota(L))\equiv \iota_*(T_xL)$ of the tangent
space $T_{\iota(x)}\R^N$ have no common nonzero vectors. This monomorphism $i$
induces the proper map $i' \equiv \pi \circ i:E \to \R^N$ of the total space
$E$ of $\xi$ into $\R^N$, sending any vector $a \in \xi_x$ to the point
$\iota(x) + i(a)$.  Suppose that this map $i'$ is transversal to $\iota(L)$
everywhere in $E \setminus L$ (where $L \subset E$ denotes the zero section of
$\xi$). Then this intersection set $ \iota(L) \cap i'(E \setminus L) $ defines
a $\Z_2$-cycle of codimension $N-l-t$ in $\iota(L) \sim L$.  \medskip

{\sc Lemma 1.} {\it The class in $H^*(L,\Z_2)$, Poincar\'e dual to this cycle
$\iota^{-1}(\iota(L) \cap i'(E \setminus L)),$ is equal to
\begin{equation}
\label{pr2} \bar w_{N-l-t}(TL \oplus \xi),
\end{equation}
the $(N-l-t)$-dimensional homogeneous component of $\bar w(TL \oplus \xi)$.

In particular, if this class (\ref{pr2}) is nontrivial, then the set $\iota(L)
\cap i'(E \setminus L)$ is not empty; moreover, the latter is true even if $i'$
is not transversal to $\iota(L)$ in $E\setminus L$.}
\medskip

{\sc Example.} If $t=0$, we get the known obstruction to the existence of an
embedding $L^l \to \R^N$, see \cite{MS}: if $\bar w_{N-l}(TL) \ne 0$, then the
{\it self-intersection set} of any such {\it embedding} is non-empty.
\medskip

{\it Proof of Lemma 1.} There is a tubular $\e$-neighborhood $U_{\e}$ of the
submanifold $\iota(L)$ in $\R^N$ such that the containing $L$ component of
$(i')^{-1}(U_{\e})$ is a tubular neighborhood of $L$ in $E$, and the
restriction of $i'$ on the latter neighborhood (which we denote by $W_\e$) is a
diffeomorphism onto its image. Let us deform the submanifold $\iota(L)$ in
$\R^N$ by a sufficiently $C^1$-small generic diffeomorphism $v$ of $\R^N$ so
that

a) $|v(y)-y|<\e/2$ for any $y \in \iota(L)$,

b) the image $v\iota(L)$ of $\iota(L)$ under this shift is a smooth manifold in
$U_\e$ transversal to the manifold $i'(W_\e)$;

c) the map $i'$ remains to be transversal to $v\iota(L)$ everywhere in $E
\setminus W_\e$.

The variety $v\iota(L) \cap i'(E)$ consists of two disjoint parts. The first,
$v\iota(L) \cap i'(E \sm W_\e)$, defines in the group $H_*(v\iota(L), \Z_2)
\simeq H_*(L, \Z_2)$ the same homology class as  $\iota(L) \cap i'(E\setminus
W_\e)$, which we wish to calculate. The second, $v\iota(L) \cap i'(W_\e)$, is
Poincar\'e dual to the first homological obstruction to the existence of a
continuous non-zero section of the quotient bundle $\R^N / (TL \oplus \xi)$
over $L$, and thus is equal to $\bar w_{N-l-t}(TL \oplus \xi)$, see \cite{MS}.

Finally, the homology class in $H_*(v\iota(L),\Z_2)$ of the sum of these two
cycles in $v\iota(L)$ is Poincar\'e dual to the restriction to $v\iota(L)$ of
the cohomology class $[E] \in H^{N-t-l}(\R^N, \Z_2)$ Poincar\'e dual to the
direct image of the fundamental cycle of the closed manifold $E$ under the
proper map $i'$. The latter class belongs to a trivial group, hence our two
cycles in $v\iota(L)$ are homologous mod 2. \quad $\Box$

\section{The general topological estimate}

Let $M$ be a $k$-dimensional manifold, and $B(M,r)$ the $r$-th configuration
space of $M$, i.e. the topological space, whose points are subsets of
cardinality $r$ in $M$. Let $\psi$ (respectively, $\tilde \psi$) be the
canonical $r$-dimensional (respectively, $(r-1)$-dimensional) vector bundle
over $B(M,r)$, whose fibre over the point $\x = \{x_1, \ldots, x_r\} \in
B(M,r)$ is the space of all functions $f: \{x_1, \ldots, x_r\} \to \R$
(respectively, of all such functions with $\sum f(x_i) = 0$). Obviously $\psi
\simeq \tilde \psi \oplus \R^1$.
\medskip

Any embedding $I:M \to \R^N$ defines the map $I^r:B(M,r) \to \R^N$, sending any
point $\x = \{x_1, \ldots, x_r\} \in B(M,r)$ to the mass center
\begin{equation}
\label{centr} \frac{1}{r}\sum I(x_i)
\end{equation}
of points $I(x_1), \ldots, I(x_r)$. The image of the tangent space
$T_{\x}B(M,r)$ under the derivative map $I^r_*$ of this map at the point $\x$
is equal to the linear span of tangent spaces $T_{I(x_i)}I(M)$, $i=1, \ldots,
r$, translated to the point (\ref{centr}).

Also there is the natural homomorphism $\chi: \tilde \psi \to T\R^N$, sending
any function $f: (x_1, \ldots, x_r) \to \R$ to the vector $\sum_{i=1}^r f(x_i)
I(x_i) \in T_{I^r(\x)}\R^N$. The corresponding subset $\pi \circ \chi(\tilde
\psi|_{\x}) \subset \R^N$ is the minimal affine subspace in $\R^N$ containing
all points $I(x_i)$.

Thus for any $\x \in B(M,r)$ we have two important subspaces  in
$T_{I^r(\x)}\R^N$:  the image of $T_\x B(M,r)$ under the derivative map $I^r_*$
and $\chi(\tilde \psi|_{\x})$.  Let $\tau(\x) \subset T_{I^r(\x)}\R^N$ be the
linear span of these two subspaces.  \medskip

{\sc Lemma 2.} {\it Suppose that the embedding $I$ is $r$-neighborly. Then for
any $\x =(x_1, \ldots, x_r) \in B(M,r)$ the affine plane $\pi(\tau(\x)) \subset
\R^N$ intersects the set $I(B(M,r))$ only at the point $I(\x)$.}
\medskip

{\it Proof.} The affine hyperplane in $\R^N$, touching $M$ at the points $x_1,
\ldots, x_r$ and participating in the definition of a $r$-neighborly embedding,
contains this plane $\pi(\tau(\x))$ but cannot contain any point of
$I^r(B(M,r))$ other than $\x$. \quad $\Box$
\medskip

{\sc Definition.} The set $\Omega(I) \subset B(M,r)$ is the set of all points
$\x$ such that $\dim \tau(\x) < kr+r-1$. For $r>1$ and $N \ge (k+1)r -1$, the
set Reg$(r) \subset C^\infty(M,\R^N)$ consists of all maps $I:M \to \R^N$ such
that the topological codimension of $\Omega(I)$ in $B(M,r)$ is greater than
$N-(k+1)r+1$ in the following exact sense: any compact
$(N-(k+1)r+1)$-dimensional submanifold in $B(M,r)$ is isotopic to one not
intersecting $\Omega(I)$. For $r=1$, set Reg$(1)= C^\infty(M,\R^N)$.  \medskip

In particular, for any such submanifold $L$, not intersecting $\Omega(I)$, the
restriction on $L$ of the bundle $\tau(\x)$ is isomorphic to $TB(M,r) \oplus
\tilde \psi$, and the restriction of the map $I$ on $L$ is an immersion into
$\R^N$ (and even an embedding if $r>1$).
\medskip

{\sc Lemma 3.} {\it For any $r>1$, the set $\Sigma \Sigma(r) \equiv
C^\infty(M,\R^N) \setminus$ Reg$(r)$ is a subset of infinite codimension in
$C^\infty (M,\R^N)$.}
\medskip

{\it Proof.} Any map $I: M \to \R^N$ defines its multijet extension $I^1_r:
B(M,r) \to J^1_r(M,\R^N)$, sending any point $(x_1, \ldots, x_r) \in B(M,r)$ to
the collection of 1-jets of $I$ at these points, see e.g. \cite{GG}. The set
$\Omega(I)$ can be described as the pre-image under this map of a certain
algebraic subset $\Sigma \subset J^1_r(M,\R^N)$, whose codimension is equal to
$N-(k+1)r+2$. If the map $I$ is of class $\Sigma \Sigma$, then this map $I^1_r$
is non-transversal to $\Sigma$ at infinitely many points, and our lemma follows
from the Thom's multijet transversality theorem, see \cite{GG}, \cite{V}. \quad
$\Box$
\medskip

{\sc Theorem 3.} {\it Suppose that $N=(k+1)r+l-1$ and there is a
$l$-dimensional compact submanifold $L \subset B(M,r)$ such that
\begin{equation}
\label{th3} \langle [L], \bar w_l (TB(M,r) \oplus \tilde \psi) \rangle \ne 0,
\end{equation}
where $[L]$ is the $\Z_2$-fundamental class of $L$. Then  there are no
$r$-neighborly embeddings $M \to \R^N$ of the class Reg$(r)$. In particular,
$\Delta(M,r) \ge (k+1)r +l.$}
\medskip

{\it Proof.} Suppose that $I$ is a $r$-neighborly embedding $M \to \R^N$ of
class Reg$(r),$ $r>1.$ Denote by $\nu$ the $(kr-l)$-dimensional vector
subbundle in the restriction of $TB(M,r)$ to $L$, orthogonal to the tangent
bundle $TL$; let be $\xi = \nu \oplus \tilde \psi_L$, so that $TL \oplus \xi =
T|_L B(M,r) \oplus \tilde \psi_L$. Since $I \in $Reg$(r)$, we can assume that
$L$ does not meet $\Omega(I)$, in particular the homomorphism $I^r_* \oplus
\chi: \xi \to T \R^N$ is a monomorphism satisfying conditions of Lemma 1. Then
conclusions of Lemmas 1 and 2 contradict to one another. Finally, in the case
$r=1$ our condition (\ref{th3}) prevents the existence of {\it any} embedding
$M \to \R^N,$ see \cite{MS} or Example after Lemma 1, and Theorem 3 is
completely proved.  \quad $\Box$

\section{Proof of main theorems}

\subsection{Proof of Proposition 1.}

Any embedding $M^k \to S^{N-1}$ obviously is a stable $1$-neighborly embedding
$M \to \R^N$. Conversely, suppose that $M$ is a 1-neighborly submanifold in
$\R^N$. The Gauss map establishes the natural homeomorphism between $S^{N-1}$
and the set of all supporting hyperplanes of $M$. For any point $x \in M$
consider the set $H(x)$ of all oriented hyperplanes, supporting $M$ only at
this point. This is a convex semialgebraic subset of the sphere $S^{N-1-k}
\subset S^{N-1}$, consisting of all oriented hyperplanes parallel to $T_xM$. If
$M$ is generic, then the set of interior points of $H(x)$ in $S^{N-1-k}$ is
non-empty, and the union of such sets forms a smooth fiber bundle over $M$ with
fibers homeomorphic to $\R^{N-1-k}$. Any smooth section of this bundle is the
desired embedding. \quad $\Box$ \medskip

{\sc Remark.} The same considerations prove that if $M$ is a (non-stably)
1-neighborly submanifold of $\R^N$, then it is {\em homeomorphic} to a subset
of $S^{N-1}$. However, the example of the curve $t \to (t, t^3, t^4)$ shows
that for non-generic embeddings the set $H(t)$ can consist of unique point, and
our fiber bundle $\{H(x) \to x \}$ can be not smooth and even not locally
trivial.

\subsection{Proof of Theorem 1!.}

First suppose that $r=2^s$. The corresponding submanifold $L(r) \subset
B(\R^k,r)$, satisfying the conditions of Theorem 3, is constructed as follows
(cf. \cite{Fuchs}, \cite{V}).

Let us fix some $\e \in (0,1/2]$. Consider in $\R^k$ a sphere of radius 1 and
mark on it some two opposite points $A_1, A_2$. Then consider two spheres of
radius $\e$ with centers in these two points and mark on any of them two
opposite points: $A_{11}, A_{12}$ on the first and $A_{21}, A_{22}$ on the
second. Consider four spheres of radius $\e^2$ with centers at all these four
points and mark on any of them two opposite points: $A_{111}, \ldots, A_{222}$,
etc. After the $s$-th step we obtain $2^s=r$ pairwise different points in
$\R^k,$ i.e. a point of the space $B(\R^k,r)$. $L(r)$ is defined as the union
of all points of the latter space, which can be obtained in this way; this is a
$(k-1)(r-1)$-dimensional smooth compact manifold.
\medskip

{\sc Proposition 2.} {\it 1. For any $k$ and $r$, the tangent bundle
$TB(\R^k,r)$ is isomorphic to the direct sum of $k$ copies of the bundle
$\psi$.

2. If $k$ is a power of 2, then the $k$th power of the total Stiefel--Whitney
class of the bundle $\tilde \psi$ (or $\psi$) over $B(\R^k,r)$ is equal to
${\bf 1} \in H^0(B(\R^k,r),\Z_2)$. In particular, $\bar w(\tilde \psi) =
(w(\tilde \psi))^{k-1}$.

3. If both $k$ and $r$ are powers of 2, then the $(k-1)(r-1)$-dimensional
homogeneous component $\bar w_{(k-1)(r-1)}(\tilde \psi) \in
H^{(k-1)(r-1)}(B(\R^k,r),\Z_2)$ of the total inverse Stiefel--Whitney class
$\bar w_*(\tilde \psi)$ is nontrivial and its value on the fundamental cycle of
the manifold $L(r)$ is equal to 1.}
\medskip

{\it Proof.} Statement 1 is obvious (and remains true if we replace $\R^k$ by
any parallelizable manifold). Statement 2 (and, moreover, similar assertion
concerning any class of the form ${\bf 1} + \{$terms of positive dimension$\}
\in H^*(B(\R^k,r),\Z_2)$) is proved in \cite{Hung}. Statement 3 is proved in
\cite{Vfil}, \S \ I.3.7 (and no doubt was known to the author of \cite{Hung}).
\quad $\Box$
\medskip

Thus the submanifold $L(r) \subset B(\R^k,r),$ $k=2^j,$ satisfies the
conditions of Theorem 3, and Theorem 1! is proved in the case when $r$ is a
power of 2.  For an arbitrary $r$ the similar submanifold $L(r)$ is constructed
as follows.  Suppose that $r = 2^{t_1} + \cdots + 2^{t_d}$, $t_1 > \cdots >
t_d,$ $d=d(r)$. Then the manifold $L(r)$ consists of all collections of $r$
points in $\R^k$, such that the collection of first $2^{t_1}$ of them belongs
to the manifold $L(2^{t_1})$, the next $2^{t_2}$ are obtained from some
collection $\x_2 \in L(2^{t_2})$ by the translation along the vector $(3,0,
\ldots,0)$, the next $2^{t_3}$ are obtained from some collection $\x_3 \in
L(2^{t_3})$ by the translation along the vector $(6,0, \ldots,0)$, etc.  In
particular, $L(r) \sim L(2^{t_1}) \times \cdots \times L(2^{t_d}).$ (If
$t_d=1$, then we set $L(t_d) = \{$the point $0\}$.) In restriction to $L(r)$,
the bundle $\tilde \psi$ is obviously isomorphic to the direct sum of $d$
bundles, induced from similar $(2^{d_i}-1)$-dimensional bundles over the
factors $L(2^{t_i})$, and the $(r-1)$-dimensional trivial bundle.  Thus the
manifold $L(r)$ for an arbitrary $r$ also satisfies conditions of Theorem 3
(with $l= (k-1)(r-d)$). \quad $\Box$

\subsection{Proof of Theorem 2!.}

Suppose that $k=2^j$. Define the submanifold $L(k,r) \subset B(\RP^k,r)$ as the
set of all unordered collections of $r$ pairwise orthogonal points in $\RP^k$
(with respect to any Euclidean metrics in $\R^{k+1}$) lying in some fixed
subspace $\RP^{k-1} \subset \RP^k$. This is a smooth $\left(kr-{r \choose
2}\right)$-dimensional manifold.

Consider also the manifold $\Lambda(k,r)$, consisting of similar ordered
collections; it is a submanifold of the space $(\RP^k)^r$ and also the space of
a $r!$-fold covering $\theta:\Lambda(k,r) \to L(k,r)$.
\medskip

{\sc Proposition 3.} {\it 1. There is a ring isomorphism
\begin{equation}
\label{ss} H^*(\Lambda(k,r),\Z_2) \simeq H^*(\RP^{k-1}, \Z_2) \otimes
H^*(\RP^{k-2}, \Z_2) \otimes \cdots \otimes H^*(\RP^{k-r}, \Z_2).
\end{equation}

2. The vector bundle over $\Lambda(k,r)$, induced by the map $\theta$ from the
bundle $\psi$ or $\tilde \psi$ on $L(k,r)$, is equivalent to the trivial one.

3. The vector bundle over $\Lambda(k,r)$, induced by the map $\theta$ from the
tangent bundle $TB(\RP^k,r)$, coincides with the restriction on $\Lambda(k,r)$
of the tangent bundle of the manifold $(\RP^k)^r$.

4. The inverse Stiefel--Whitney class $\bar w_* \equiv (w_*)^{-1}$ of the
latter tangent bundle satisfies the inequality}
$$\langle [\Lambda(k,r)], \bar w_{kr-{r \choose 2}}(T(\RP^k)^r)
\rangle \ne 0.$$

{\it Proof.} Statement 1 is a standard exercise on homology of fiber bundles,
see e.g. \cite{Borel}. Namely, consider the fiber bundle $p: \Lambda(k,r) \to
\Lambda(k,r-1)$, sending any ordered collection $(x_1, \ldots, x_r)$ to $(x_1,
\ldots, x_{r-1})$. Its fiber $\F$ is equal to $\RP^{k-r}$, hence the
fundamental group of the base acts trivially on $H^*(\F,\Z_2)$, and the term
$E_2^{p,q}$ of the $\Z_2$-spectral sequence of this bundle is naturally
isomorphic to $H^p(\Lambda(k,r-1),\Z_2)\otimes H^q(\F,\Z_2)$. All further
differentials $d^i$, $i \ge 2$, of the spectral sequence act trivially on all
elements of the column $E^{0,*} \simeq H^*(\F,\Z_2)$, because the embedding
homomorphism $H^*(\Lambda(k,r), \Z_2) \to H^*(\F,\Z_2) $ is epimorphic: indeed,
its composition with the map $H^*(\RP^k, \Z_2) \to H^*(\Lambda(k,r), \Z_2),$
defined by the projection of the space $\Lambda(k,r) \subset (\RP^k)^r$ onto
the $r$-th copy of $\RP^k$, is just the epimorphism $H^*(\RP^k, \Z_2) \to
H^*(\RP^{k-r}, \Z_2)$, induced by the embedding. Since the spectral sequence is
multiplicative, it degenerates at the term $E_2$.

Statements 2 and 3 of Proposition 3 are obvious, and statement 4 follows from
the induction conjecture for $\Lambda(k,r-1)$ and the fact that for $k=2^j$ the
inverse Stiefel--Whitney class $\bar w_*(T\RP^k)$ is equal to
$1+\alpha+\alpha^2 + \cdots + \alpha^{k-1}$, where $\alpha$ is the
multiplicative generator of the ring $H^*(\RP^k, \Z_2)$, see e.g. \cite{MS}.
\quad $\Box$
\medskip

Using the functoriality of Stiefel--Whitney classes and the Whitney
multiplication formula for these classes of a Whitney sum of bundles, we obtain
from this proposition that $\bar w_{kr-{r \choose 2}}(TB(\RP^k,r) \oplus \tilde
\psi)$ is a non-trivial element of $H^{kr-{r \choose 2}}(L(k,r),\Z_2)$. Theorem
2! is now reduced to Theorem 3.  \quad $\Box$
\medskip

{\sc Remark.} The number $\Delta(M,k)$ is a measure of the ``topological
complexity'' of the configuration space $B(M,r)$. This complexity appears from
two issues: the obvious free action of the symmetric group $S(r)$ and the
topological complexity of the manifold $M$ itself. In Theorem 1! we essentially
exploit only the first issue, and in Theorem 2! only the second.  The
simultaneous consideration of these two interacting components should give us
more precise estimates of $\Delta(M,r)$.

\end{document}